\documentclass[12pt]{amsart}
\usepackage{amsfonts}
\usepackage{amssymb}
\theoremstyle{definition}

\theoremstyle{remark}

\newcommand{\ds}{\displaystyle}

\begin{document}

\centerline{\large\bf ALMOST HERMITIAN MANIFOLDS WITH VANISHING}
\centerline{\large\bf GENERALIZED BOCHNER CURVATURE TENSOR
\footnote{\it SERDICA Bulgaricae mathematicae publicationes. Vol. 9, 1983, p. 96-101}}

\vspace{0.2in}
\centerline{\large OGNIAN T. KASSABOV}

\vspace{0.5in}
{\sl We deal with the generalized Bochner curvature tensor and the Bochner curvature tensor
introduced respectively in [1] and [11]. In section 2 we prove, that if an almost Hermitian
manifold is a product of two almost Hermitian manifolds $M_1$ and $M_2$, then $M_1$ (resp.
$M_2$) is of pointwise constant holomorphic sectional curvature $\mu$ (resp. $-\mu$). In
section 3 we obtain a classification theorem for non K\"ahler nearly K\"ahler manifolds
with vanishing Bochner curvature tensor and a classification theorem for nearly K\"ahler manifolds
with vanishing Bochner curvature tensor and constant scalar curvature.}

\vspace{0.2in}
{\bf 1. Preliminaries.} Let $M$ be a $2m$ - dimensional almost Hermitian manifolds with
Riemannian metric $g$ and almost complex structure $J$ and let $\nabla$ be the covariant
differentiation on $M$. The curvature tensor $R$ is defined by
$$
    R(X,Y,Z,U) = g(\nabla_X\nabla_Y Z - \nabla_Y\nabla_X Z - \nabla_{[X,Y]} Z,U)
$$
for $X,\, Y,\, Z,\, U \in \mathfrak X(M)$. From the curvature tensor \ $R$ \ one may construct a tensor
$R^*$ [1] by
$$
    R^*(X,Y,Z,U) =\frac{3}{16} \{ R(X,Y,Z,U) + R(X,Y,JZ,JU)
$$
$$
    + R(JX,JY,Z,U) + R(JX,JY,JZ,JU) \}
$$
$$
    +\frac{1}{16} \{ R(JX,JZ,Y,U) - R(JY,JZ,X,U) + R(X,Z,JY,JU) - R(Y,Z,JX,JU)
$$
$$
    + R(Y,JZ,JX,U) - R(X,JZ,JY,U) +R(JY,Z,X,JU) - R(JX,Z,Y,JU) \} .
$$
The tensor $R^*$ has the following properties:

1) $R^*(X,Y,Z,U) = - R^*(Y,X,Z,U)$,

2) $R^*(X,Y,Z,U) = -R^*(X,Y,U,Z)$,

3) $\underset{X,Y,Z}\sigma R^*(X,Y,Z,U) = 0 $,

\noindent
where \ $\sigma$ \ denotes the cyclic sum and \ $R^*$ \ is the only tensor with these
properties for which
$$
    R^*(X,Y,Z,U) = R^*(JX,JY,Z,U) \ , \ \ R^*(X,JX,JX,X) = R(X,JX,JX,X) .
$$

Let \ $\{ E_i;\, i=1,...,2m\}$ \ be an orthonormal local frame field. The Ricci tensor
\ $S$ \ and the scalar curvature \ $\tau(R)$ \ of \ $M$ \ are defined by
$$
    S(X,Y) = \sum_{i=1}^{2m} R(X,E_i,E_i,Y) \ , \ \ \tau(R) = \sum_{i=1}^{2m} S(E_i,E_i).
$$
Analogously one denotes
$$
    S^*(X,Y) = \sum_{i=1}^{2m} R^*(X,E_i,E_i,Y) \ , \ \ \tau^*(R) = \sum_{i=1}^{2m} S^*(E_i,E_i)
$$
and it is easy to see that \ $S^*(X,Y)=S^*(JX,JY) = S^*(Y,X)$.

The generalized Bochner curvature tensor  $B^*$  for  $M$  is defined by [1]
$$
    B^*=R^* - \frac 1{2(m+2)} (\varphi + \psi)(S^*) +
        \frac{\tau^*(R)}{4(m+1)(m+2)} (\pi_1+\pi_2),
$$
where
$$
    \varphi(Q)(X,Y,Z,U) = g(X,U)Q(Y,Z) - g(X,Z)Q(Y,U)
$$
$$
    + g(Y,Z)Q(X,U) - g(Y,U)Q(X,Z)\, ,
$$
$$
    \psi(Q)(X,Y,Z,U) = g(X,JU)Q(Y,JZ) - g(X,JZ)Q(Y,JU)
$$
$$
    -2g(X,JY)Q(Z,JU) + g(Y,JZ)Q(X,JU) - g(Y,JU)Q(X,JZ)-2g(Z,JU)Q(X,JY)\, ,
$$
$$
    \pi_1(X,Y,Z,U) = g(X,U)g(Y,Z) - g(X,Z)g(Y,U)\, ,
$$
$$
    \pi_2(Q)(X,Y,Z,U) = g(X,JU)Q(Y,JZ) - g(X,JZ)Q(Y,JU) -2g(X,JY)Q(Z,JU)\, ,
$$
for an arbitrary tensor  $Q$  of type \ (0,2).

An almost Hermitian manifold \, $M$ \, which satisfies \, $(\nabla_X)X=0 $ \, for all \
$X \in \mathfrak X(M)$ \ is called a nearly K\"ahler manifold. It is well known that for
such a manifold \ $R(X,Y,Z,U)=R(JX,JY,JZ,JU)$ \ holds good for all
$X,\, Y,\, Z,\, U \in {\mathfrak X}(M)$ [3]. In general a manifold which satisfies this
identity is said to be an RK-manifold. For an RK-manifold \, $S(X,Y)=S(JX,JY)$ \, holds. The
Bochner curvatute tensor for an RK-manifold of dimension $2m \ge 6$ is defined [12] by
$$
    B=R-\frac 1{8(m+2)}(\varphi + \psi)(S+3S')-\frac 1{8(m-2)}(3\varphi - \psi)(S-S')
$$
$$
    +\frac{\tau(R)+3\tau'(R)}{16(m+1)(m+2)}(\pi_1+\pi_2)
    +\frac{\tau(R)-\tau'(R)}{16(m-1)(m-2)}(3\pi_1-\pi_2),
$$
where
$$
    S'(X,Y) = \sum_{i=1}^{2m} R(X,E_i,JE_i,JY) \ , \ \ \tau'(R) = \sum_{i=1}^{2m} S'(E_i,E_i).
$$
It is easy to check that for an RK-manifold \, $4S^*=S+3S'$.

Now let $M$ be a nearly K\"ahler manifold and $X,\, Y,\, Z,\, U,\, V \in {\mathfrak X}(M)$. We
shall use the following formulas (see [3; 6; 12]):
$$
    R(X,Y,Z,U) - R(X,Y,JZ,JU) = -g((\nabla_XJ)Y,(\nabla_ZJ)U),   \leqno (1.1)
$$
$$
    2g((\nabla_X(\nabla_Y J))Z,U)= \underset{Y,U,Z}{\sigma} R(X,JY,U,Z),   \leqno (1.2)
$$
$$
    2(\nabla_X(S-S'))(Y,Z)=(S-S')((\nabla_XJ)Y,JZ)+(S-S')(JY,(\nabla_XJ)Z), \leqno (1.3)
$$
$$
    X(\tau(R)-\tau(R'))=0 ,        \leqno (1.4)
$$
$$
    \sum_{i,j=1}^{2m}(S-S')(E_i,E_j)(S-5S')(E_i,E_j)=0.        \leqno (1.5)
$$
From the second Bianchi identity
$$
    \underset{X,Y,Z}{\sigma} (\nabla_X R)(Y,Z,U,V)=0
$$
it follows
$$
    \sum_{i=1}^{2m} (\nabla_{E_i}R)(X,Y,Z,E_i)=(\nabla_X S)(Y,Z)-(\nabla_Y S)(X,Z) , \leqno (1.6)
$$
$$
    \sum_{i=1}^{2m}(\nabla_{E_i}S)(X,E_i)=\frac 12 X(\tau(R)).        \leqno (1.7)
$$

\vspace{0.1in}
{\bf 2. Product of almost Hermitian manifolds and the generalized Bochner curvature tensor.}

T h e o r e m 2.1. {\it Let an almost Hermitian manifold $M$ be a product
$M_1 \times M_2$, where $M_1$ and $M_2$ are almost Hermitian manifolds. Then $M$ has
vanishing generalized Bochner curvature tensor if and only if $M_1$ (resp. $M_2$)
is of pointwise constant holomorphic sectional curvature $\mu$ (resp. $-\mu$).}

P r o o f. If $B^* = 0$ we find
$$
    \begin{array}{c}
        \ds R^*(X,Y,Z,U) = \frac 1{2(m+2)} \{g(X,U)S^*(Y,Z) - g(X,Z)S^*(Y,U) \\
                            + g(Y,Z)S^*(X,U) - g(Y,U)S^*(X,Z)+ g(X,JU)S^*(Y,JZ)        \\
                - g(X,JZ)S^*(Y,JU) + g(Y,JZ)S^*(X,JU) - g(Y,JU)S^*(X,JZ)   \\
                -2g(X,JY)S^*(Z,JU)-2g(Z,JU)S^*(X,JY) \}        \\
    \ds            - \frac{\tau^*(R)}{4(m+1)(m+2)} \{g(X,U)g(Y,Z) - g(X,Z)g(Y,U)+ g(X,JU)g(Y,JZ) \\
            -g(X,JZ)g(Y,JU)-2g(X,JY)g(Z,JU)\}.
  \end{array}        \leqno (2.1)
$$
We denote by $R_1$ the curvature tensor of $M_1$. Analogously we have $R_1^*$, $S_1^*$,
$\tau^*(R_1)$. Note that $R=R_1$, $R^*=R_1^*$ and $S^*=S_1^*$ on $M_1$. Let
$X \in {\mathfrak X}(M)$ and $\{ E_i;\ i=1,...,2k \}$ be an orthonormal local frame field on $M_1$.
In (2.1) we put $U=X$, $Y=Z=E_i$ and adding for $i=1,...,2k$ we obtain
$$
    S_1^*(X,X) = \left\{ \frac{\tau^*(R)}{2(m-k)} - \frac{(k+1)\tau^*(R)}{2(m+1)(m-k)} \right\}
                 g(X,X) . \leqno (2.2)
$$
From (2.1) it follows for a unit field $X \in {\mathfrak X}(M_1)$
$$
    R_1(X,JX,JX,X)= \frac 4{m+2}S_1^*(X,X) - \frac{\tau^*(R)}{(m+1)(m+2)}.   \leqno (2.3)
$$
Because of (2.2) and (2.3) $M_1$ is of pointwise constant holomorphic sectional curvature
$\mu$. Hence
$$
    S_1^*(X,X) = \frac{k+1}{2}\mu g(X,X)  \leqno (2.4)
$$
for $X \in {\mathfrak X}(M_1)$ and
$$
    \tau^*(R_1)=k(k+1)\mu .    \leqno (2.5)
$$
Analogously $M_2$ is of pointwise constant holomorphic sectional curvature $\mu'$ and
$$
    \tau^*(R_2)=(m-k)(m-k+1)\mu' .     \leqno (2.6)
$$
Using (2.3)-(2.6) and $\tau^*(R) = \tau^*(R_1)+\tau^*(R_2)$ we find $\mu'=-\mu$.

The converse is a simple calculation by using the fact, that $M$ is of pointwise
constant holomorphic sectional curvature $\mu$ if and only if
$R^* = \frac {\mu}4(\pi_1 + \pi_2)$, see [1].

C o r o l l a r y \, 2.2. {\it Let an almost Hermitian manifold $M$ be a product of
more than two almost Hermitian manifolds. Then $M$ has vanishing generalized Bochner
curvature tensor if and only if it is of zero holomorphic sectional curvature.}

For the case of a K\"ahler manifold $M$ in Theorem 2.1 and Corollary 2.2 see [10].

\vspace{0.1in}
{\bf 3. Nearly K\"ahler manifolds with vanishing Bochner curvature tensor.}

L e m m a. {\it \, Let \, $M$ \, be a \, $2m$-dimensional nearly K\"ahler manifold, \, $m>2$. If
the Bochner curvature tensor of \, $M$ vanishis, then the tensor \, $S-S'$ \, is
parallel.}

P r o o f. Let \ $\{ E_i;\ i=1,...,2m \}$ be a local orthonormal frame field. From (1.1) we find
$$
    (S-S')(X,Y)=\sum_{i=1}^{2m} g((\nabla_X J)E_i,(\nabla_Y J)E_i) ,
$$
which implies
$$
    (\nabla_X(S-S'))(Y,Y)= 2\sum_{i=1}^{2m} g((\nabla_X(\nabla_Y J))E_i,(\nabla_Y J)E_i) .
$$
Hence, using (1.2) we obtain
$$
    \begin{array}{c}
        \ds (\nabla_X(S-S'))(Y,Y)= \sum_{i=1}^{2m} \{ R(X,JY,(\nabla_YJ)E_i,E_i) \\
        +R(X,J(\nabla_YJ)E_i,E_i,Y) + R(X,JE_i,Y,(\nabla_Y J)E_i) \}.
    \end{array}   \leqno (3.1)
$$
From $B=0$, (3.1) and (1.3) we derive $\nabla (S-S') =0$.

T h e o r e m \, 3.1. \, {\it Let $M$ be a $2m$-dimensional non K\"ahler nearly K\"ahler
manifold, $m>2$. If $M$ has vanishing Bochner curvature tensor, it is locally
isometric to one of:

a) the sixth sphere $ {\bf S}^6$;

b) $ {\bf CD}^1(-c)\times {\bf S}^6(c)$ \, where  $ {\bf CD}^1(-c)$ \, (resp. \, $ {\bf S}^6(c)$) \,
is the one-dimensional complex hyperbolic space of constant sectional curvature $-c$ (resp. the
sixth sphere of constant sectional curvature $c$).}

P r o o f. According to the lemma the tensor $S-S'$ is parallel. Let $M$ be locally a product \,
$M_1(\lambda_1)\times ... \times M_k(\lambda_k)$, where \, $S-S'= \lambda_ig$ \, on \, $M_i(\lambda_i)$ \,
and \, $\lambda_i \ne \lambda_j$ \, for  $i \ne j$. As it's easy to see $M_i(\lambda_i)$ is a
nearly K\"ahler manifold for $i=1,...,k$. From $B=0$ it follows $B^*=0$. So Theorem 2.1, Corollary 2.2
and [9] imply that if $k>1$ \, $M$ is either of zero holomorphic sectional curvature, or locally a
product of two nearly K\"ahler manifolds of constant holomorphic sectional curvature $-\mu$ and $\mu$,
respectively, $\mu>0$. Since $M$ is non K\"ahler, the former is impossible and the latter occurs only
when $M$ is locally isometric to \, $ {\bf CD}^{m-3}(-c)\times {\bf S}^6(c)$ [5]. On the other hand
\, $ {\bf CD}^{m-3}(-c)\times {\bf S}^6(c)$ has nonvanishing Bochner curvature tensor, if $m>4$.
Indeed it doesn't satisfy the condition $R(x,y,z,u)=0$ for all mutually orthogonal $x$, $y$, $z$, $u$ spanning a
4-dimensional antiholomorphic plane.

Let $k=1$. Then
$$
    S-S'=\frac{1}{2m} \{ \tau(R)-\tau'(R) \}g.     \leqno (3.2)
$$
Now (3.2) and (1.5) imply
$$
    (\tau(R)-\tau'(R))(\tau(R)-5\tau'(R))=0.
$$
If $\tau(R)(p)-\tau'(R)(p)=0$ for some $p \in M$ then $\tau(R)-\tau'(R)=0$ because of (1.4).
Then from (3.2) and $B=0$ we obtain
$$
    R = \frac 1{2(m+2)}(\varphi +\psi)(S) - \frac{\tau(R)}{4(m+1)(m+2)}(\pi_1+\pi_2).
$$
Hence \, $R(X,Y,Z,U) =R(X,Y,JZ,JU)$ \, holds for all \, $X,\, Y,\, Z,\, U \in {\mathfrak X}(M)$
and so $M$ is a K\"ahler manifold [4], which is a contradiction. Consequently
$$
    \tau(R)-5\tau'(R)=0.    \leqno (3.3)
$$
Now (1.4) implies that $\tau(R)$ and $\tau'(R)$ are global constants. From $B=0$, (3.2)
and (3.3) we find
$$
    R=\frac{1}{2(m+2)} (\varphi+\psi)(S) - \frac{(4m+3)\tau(R)}{10m(m+1)(m+2)} (\pi_1 +\pi_2)
                                         + \frac{\tau(R)}{20m(m-1)} (3\pi_1 -\pi_2).  \leqno(3.4)
$$
Using (1.6), (1.7), (3.4) and $X(\tau(R))=0$ we obtain
$$
    (2m+3)\{ (\nabla_X S)(Y,Y) - (\nabla_Y S)(X,Y) \} =
       3(\nabla_{JY} S)(X,JY) - 3 S((\nabla_X J)Y,JY). \leqno (3.5)
$$
In particular \, $(\nabla_X S)(X,X)=0$ \, and hence
$$
    (\nabla_X S)(Y,Y) + 2(\nabla_Y S)(X,Y)=0.    \leqno (3.6)
$$
From (3.5) and (3.6) we  derive \, $(\nabla_X S)(Y,Y)=0$ \, which implies \, $\nabla S=0$.
Since $M$ is not locally isometric to \, $ {\bf CD}^1(-c)\times {\bf S}^6(c)$, it
follows that \, $ S=\tau(R)g/(2m)$. Now (3.4) takes the form
$$
    R = \frac{5m+1}{20m(m^2-1)}\tau(R)\pi_1 + \frac{m-3}{20m(m^2-1)}\tau(R)\pi_1.
$$
Consequently $M$ is of pointwise constant holomorphic sectional curvature. Since $M$ is
not K\"ahler, it is locally isometric to ${\bf S}^6$ [5].

From [8] and Theorem 3.1 we obtain

T h e o r e m \, 3.2. \, {\it Let $M$ be a $2m$- dimensional nearly K\"ahler manifold, $m>2$,
with vanishing Bochner curvature tensor and constant scalar curvature. Then $M$ is locally
isometric to one of the following:

a) the complex Euclidian space ${\bf CE}^m$;

b) the complex hyperbolic space ${\bf CD}^m$;

c) the complex projective space ${\bf CP}^m$;

d) the sixth sphere $ {\bf S}^6$;

e) the product $ {\bf CD}^1(-c)\times {\bf S}^6(c)$;

f) the product $ {\bf CD}^{m_1}(-c)\times {\bf CP}^{m_2}(c)$, $m_1+m_2=m$.}

We note that Theorems 3.1 and 3.2 can be proved also by using Theorem 2.1,
Corollary 2.2 and [6, Theorem 4.11].

Since for an almost Hermitian manifold the condition of constant antiholomorphic
sectional curvature implies $B=0$, see [7], we obtain

C o r o l l a r y \, 3.3 \, [2]. \, {\it Let $M$ be a $2m$-dimensional nearly K\"ahler
manifold, $m>2$. If $M$ is of pointwise constant antiholomorphic sectional
curvature $\nu$, it is locally isometric to one of the following:

a) the complex Euclidian space ${\bf CE}^m$;

b) the complex hyperbolic space ${\bf CD}^m(4\nu)$;

c) the complex projective space ${\bf CP}^m(4\nu)$;

d) the sixth sphere $ {\bf S}^6(\nu)$.}

\vspace{0.4in}
\centerline{\large R E F E R E N C E S}

\vspace{0.2in}
\noindent
1. G. G a n \v{c} e v. Characteristic of some classes of almost Hermitian manifolds.
{\it Serdica},

{\bf 4}, 1978, 19-23.

\noindent
2. G. T. G a n \v{c} e v, \, O. T. K a s s a b o v. Nearly K\"ahler manifolds of constant antiho-

lomorphic sectional curvature. {\it C. R. Acad. bulg. Sci.}, {\bf 35}, 1982, 145-147.

\noindent
3. A. G r a y. Nearly K\"ahler manifolds. {\it J. Diff. Geom.}, {\bf 4}, 1970, 283-309.

\noindent
4. A. G r a y. Vector cross products on manifolds. {\it Trans. Amer. Math. Soc.}, {\bf 141}, 1969,

465-504.

\noindent
5. A. G r a y. Classification des vari\'et\'es approximativement k\"ahleriennes de courbure

sectionnelle holomorphe constante. {\it C. R. Acad. Sci. Paris, S\'er. A}, {\bf 279}, 1974, 797-800.

\noindent
6. A. G r a y. The structure of nearly K\"ahler manifolds. {\it Math. Ann.}, {\bf 223}, 1976, 233-248.

\noindent
7. O. K a s s a b o v. On the Bochner curvature tensor in an almost Hermitian manifold.

{\it Serdica}, {\bf 9}, 1983, , 168-171.

\noindent
8. M. M a t s u m o t o, \, S. T a n n o. K\"ahlerian spaces with parallel or vanishing

Bochner curvature tensor. {\it Tensor} (N.S.), {\bf 27}, 1973, 291-294.

\noindent
9. A. M. N a v e i r a, \, L. M. H e r v e l l a. Schur's theorem for nearly K\"ahler manifolds.

{\it Proc. Amer. Math. Soc.}, {\bf 49}, 1975, 421-425.

\noindent
10. S. T a c h i b a n a, \, R. L i u. Notes on K\"ahlerian metrics with vanishing Bochner

\ curvature tensor. {\it K\"odai Math. Sem. Rep.}, {\bf 22}, 1970, 313-321.

\noindent
11. F. T r i c e r r i, \, L. V a n h e c k e. Curvature tensors on almost Hermitian manifolds.

\ {\it Trans. Amer. Math. Soc.}, {\bf 267}, 1981, 365-398.

\noindent
12. Y. W a t a n a b e, \, K. T a k a m a t s u. \, On a $K$-space of constant holomorphic

\ sectional curvature. {\it K\"odai Math. Sem. Rep.}, {\bf 25}, 1972, 351-354.

\vspace {0.5cm}
\noindent
{\it Center for mathematics and mechanics \ \ \ \ \ \ \ \ \ \ \ \ \ \ \ \ \ \ \ \ \ \ \ \ \ \ \ \ \ \ \ \ \ \
Received 10.4.1981

\noindent
1090 Sofia   \ \ \ \ \ \ \ \ \ \ \ \ \ \ \ \ \  P. O. Box 373}

\end{document}